\documentclass[12pt,draft]{amsart}
\usepackage{amsmath,amsthm,latexsym,amscd,amsbsy,amssymb}

\chardef\bslash=`\\ 

\makeatletter
\def\verbatim{\interlinepenalty\@M \@verbatim
  \leftskip\@totalleftmargin\advance\leftskip2pc
  \frenchspacing\@vobeyspaces \@xverbatim}
\makeatother
\newtheorem{thm}{Theorem}[section]
\newtheorem{cor}[thm]{Corollary}
\newtheorem{lem}[thm]{Lemma}
\newtheorem{pro}[thm]{Proposition}

\newtheorem*{A}{Theorem A}
\newtheorem*{B}{Theorem B}
\newtheorem*{C}{Theorem C}
\newtheorem*{D}{Theorem D}
\newtheorem*{E}{Theorem E}

\theoremstyle{definition}

\theoremstyle{remark}

\numberwithin{equation}{section}

\font\f=msbm10

\begin{document}


\title[Torus groups]
{Topological characterization of torus groups}
\author{Alex Chigogidze}
\address{Department of Mathematics and Statistics,
University of Saskatche\-wan,
McLean Hall, 106 Wiggins Road, Saskatoon, SK, S7N 5E6,
Canada}
\email{chigogid@math.usask.ca}
\thanks{Author was partially supported by NSERC research grant.}

\dedicatory{To the memory of Professor Ball}

\begin{abstract}{Topological characterization of torus groups  is given.}
\end{abstract}

\maketitle
\markboth{A.~Chigogidze}{Torus groups}

\section{Introduction}\label{S:intro}
The problem of topological characterization of torus groups (i.e., powers of the circle group {\f T}) has been of considerable importance in the theory of compact Abelian groups. The following statement collects most of the known information.

\begin{A}\label{T:A}
Let $G$ be a {\em connected} compact Abelian group. Then the following conditions are equivalent:
\begin{itemize}
\item[(1)]
$G$ is arcwise connected.
\item[(2)]
$G$ is locally arcwise connected.
\item[(3)]
The map $\exp_{G} \colon L(G) \to G$ is surjective.
\item[(4)]
The map $\exp_{G} \colon L(G) \to G$ is open.
\item[(5)]
$\operatorname{Ext}(\widehat{G},\text{\f Z}) = 0$, i.e., the character group $\widehat{G}$ of $G$ is a Whitehead group.
\item[(6)]
$\operatorname{Ext}(\widehat{G},\text{\f Z}) = 0$ and each pure subgroup of $\widehat{G}$ of finite rank splits.
\end{itemize}
In addition, if $G$ is {\em metrizable}, then the above conditions are equivalent to each of the following:
\begin{itemize}
\item[(7)]
$\widehat{G}$ is free.
\item[(8)]
$\widehat{G}$ is $\aleph_{1}$-free, i.e., every countable subgroup of $\widehat{G}$ is free.
\item[(9)]
$\widehat{G}$ is projective.
\item[(10)]
$G$ is injective (in the category of compact Abelian groups).
\item[(11)]
$G$ is locally connected.
\item[(12)]
$G$ is a torus.
\end{itemize}
\end{A}

In the nonmetrizable case situation is quite different. For instance, the equivalence $(8) \Longleftrightarrow (11)$ remains true, whereas $(11) \Longleftrightarrow (12)$ does not. The equivalence $(5) \Longleftrightarrow (7)$ is undecidable in ZFC \cite{shelah}. In topological terms this means that:
\begin{itemize}
\item[(a)]
There are (necessarily nonmetrizable) connected and locally connected compact Abelian groups which are not tori.
\item[(b)]
It is undecidable within ZFC that the arcwise connectedness of a nonmetrizable compact Abelian group forces the group to be isomorphic to a torus.
\end{itemize}

Below we show (Theorem E) that a simple topological property gives a complete characterization of torus groups. Statements of this sort are not unusual. Here are some examples (see also \cite{keesl}).
\begin{B}[Anderson-Kadec \cite{anderson},\cite{kadec} $\tau = \omega$; \cite{chi91} $\tau > \omega$]\label{T:B}
The following conditions are equivalent for an infinite-dimensional locally convex linear topological space  $E$ of weight $\tau \geq \omega$:
\begin{itemize}
\item[(a)]
$E$ is topologically equivalent to $R^{\tau}$ (the $\tau$-th power of the real line).
\item[(b)]
$E$ is an $AE(0)$-space.
\end{itemize}
\end{B}

\begin{C}[Pontryagin-Haydon, \cite{pontr}, \cite{haydon}]\label{T:C}
Every compact group is an $AE(0)$-space.
\end{C}

\begin{D}[\cite{bellchi}]\label{T:D}
The following conditions are equivalent for a zero-dimensional topological group $G$:
\begin{itemize}
\item[(a)]
$G$ is topologically equivalent to the product $(\text{\f Z}_{2})^{\tau} \times \text{\f Z}^{\kappa}$.
\item[(b)]
$G$ is an $AE(0)$-space.
\end{itemize}
\end{D}

Our goal in these notes is to prove the following statement.
\begin{E}\label{T:E}
The following conditions are equivalent for a compact Abelian group $G$:
\begin{itemize}
\item[(a)]
$G$ is a torus group (both in topological and algebraic senses).
\item[(b)]
$G$ is an $AE(1)$-compactum.
\item[(c)]
The map $\exp_{G} \colon L(G) \to G$ is $0$-soft.
\end{itemize}
\end{E}

There are various equivalent versions of definitions of the concepts of $AE(n)$-spaces and $n$-soft maps. In this particular situation categorical terminology seems to be more appropriate. Suppose that $\boldsymbol{\mathcal {B}}$ is a full subcategory of a category $\boldsymbol{\mathcal {A}}$. An object $a \in \operatorname{Ob}(\boldsymbol{\mathcal {A}})$ is called $\boldsymbol{\mathcal {B}}$-injective, if for any monomorphism $i \in \operatorname{Mor}_{\boldsymbol{\mathcal {B}}}(b,c)$ of the category $\boldsymbol{\mathcal {B}}$ and for any morphism $f \in \operatorname{Mor}_{\boldsymbol{\mathcal {A}}}(b,a)$ there exists a morphism $\tilde{f} \in \operatorname{Mor}_{\boldsymbol{\mathcal {A}}}(c,a)$ such that $f = \tilde{f}i$.

Let now ${\boldsymbol{\mathcal {TYCH}}}$ denote the category of Tychonov spaces and their continuous maps and ${\boldsymbol{\mathcal {TYCH}}_{n}}$ be the full category consisting of at most $n$-dimensional (in the sense of $\dim$) spaces. Then $AE(n)$-spaces are precisely ${\boldsymbol{\mathcal {TYCH}}_{n}}$-injective objects of ${\boldsymbol{\mathcal {TYCH}}}$. $n$-soft maps can be defined similarly: these are precisely $\boldsymbol{Mor}(\boldsymbol{\mathcal {TYCH}}_{n})$-injective objects of the category $\boldsymbol{Mor}(\boldsymbol{\mathcal {TYCH}})$. The class of metrizable $AE(0)$-spaces coincides with the class of separable completely metrizable spaces and in positive dimensions the strictly decreasing sequences
\[ AE(1) \supset AE(2) \supset \dots \supset AE(n) \supset AE(n+1) \supset \dots \]
and 
\[ LC^{0}\cap C^{0} \supset LC^{1}\cap C^{1} \dots \dots LC^{n-1}\cap C^{n-1} \supset LC^{n}\cap C^{n} \supset \dots \]
are identical (i.e., $AE(n) = LC^{n-1}\cap C^{n-1}$) in  the presence of metrizability and differ (i.e. $AE(n) \subsetneq LC^{n-1} \cap C^{n-1}$) in general. Important examples of  nonmetrizable $AE(0)$-spaces are uncountable powers of the real line. This fact is exploited below when considering the above mentioned space $L(G) \cong \text{\f R}^{\operatorname{rank}(\widehat{G})}$. Every $0$-soft map is open and surjective and for maps between completely metrizable separable spaces the converse is also true. This fact is also exploited when considering the map $\exp_{G} \colon L(G) \to G$.

Since, by duality, any statement in the category of compact Abelian groups generates an equivalent statement in the category of Abelian groups (see Theorem A), it might be worth examining Theorem E from this point of view. By the equivalence $(7) \Longleftrightarrow (12)$ in Theorem A, the statement corresponding to Theorem E in the category of Abelian groups must provide a characterization of free groups. Since $\dim G = \operatorname{rank}(\widehat{G})$ for a finite dimensional compact Abelian group $G$ and since the concept of $AE(1)$ is defined in terms $1$-dimensional testing compacta, this would lead us to the concept of $1$-projectivness (restrict, in the standard definition of projectivness, domains of testing epimorphisms by groups of rank $1$) for Abelian groups and to a hope of subsequent characterization of free Abelian groups in these terms. The so obtained version of $AE(1)$ will differ from the original topological concept of $AE(1)$ and will not be a replacement of the latter in Theorem E. Reason is simple - every $1$-dimensional connected compact (Abelian) group is metrizable and, consequently, the changed version of $AE(1)$ would essentially coincide with the arcwise connectedness. This shows that Theorem E provides a purely {\em topological} characterization of tori and leaves the problem of characterizing free uncountable Abelian groups until Section \ref{S:free}.

\section{Preliminaries}\label{S:prelim}
All spaces below are Tychonov (i.e., completely regular and Hausdorff) and maps are continuous. {\f R}, {\f Z} and $\text{\f T} = \text{\f R}/\text{\f Z}$ stand for the additive group of real numbers, the integers and the circle group respectively. We assume the familiarity with the general theory of $ANE(n)$-spaces and $n$-soft maps as well as with the spectral machinery based on the \v{S}\v{c}epin's Spectral Theorem \cite{shch} (a comprehensive discussion of these topics including a noncompact version of the Spectral Theorem can be found in \cite{book}). The space $L(G)$, defined \cite{chm} by letting $L(G) \stackrel{\operatorname{def}}{=} \operatorname{Hom}(\text{\f R},G)$, is a real topological linear space (with the topology of uniform convergence on compact sets) algebraically and topologically isomorphic to $\text{\f R}^{\operatorname{rank}(\widehat{G})}$. The continuous homomorphism $\exp_{G} \colon L(G) \to G$ is defined naturally: $\exp_{G}(f) = f(1)$, $f \in L(G)$. If $p \colon G \to T$ is a continuous homomorphism of compact Abelian groups, then the map $L(p) \colon L(G) \to L(T)$, defined by letting $L(p)(f) = pf$, $f \in L(G)$, is a continuous homomorphism.

As was noted above $X$ is an $AE(1)$-compactum if any map $g_{0} \colon Z_{0} \to X$, defined on a closed subspace $Z_{0}$ of an at most $1$-dimensional compactum $Z$, admits an extension $g \colon Z \to X$. Note that the arcwise connectedness of $X$ is obtained by restricting the previous definition to the one particular case: $Z_{0} = \{ 0,1\}$ and $Z = [0,1]$. $1$-soft maps are defined similarly. A map $f \colon X \to Y$ between compacta is $1$-soft if for any maps $g_{0} \colon Z_{0} \to X$ and $h \colon Z \to Y$ such that $fg_{0} = h$, where $Z_{0}$ is a closed subspace of an at most $1$-dimensional compactum $Z$, there exists a map $g \colon Z \to X$ such that $g_{0} = g|Z_{0}$ and $h = fg$. 

When a well-ordered inverse spectrum ${\mathcal S}_{G} = \{ G_{\alpha}, p_{\alpha}^{\alpha +1},\tau \}$ is used to investigate properties of its limit, a crucial information is contained not only in spaces $G_{\alpha}$ but in short projections $p_{\alpha}^{\alpha +1} \colon G_{\alpha +1} \to G_{\alpha}$ as well. In other words, {\em in order to obtain information about the limit object one needs information about morphisms of the spectrum representing the given object}. This shows how Lemma \ref{L:parametric} is related to Proposition \ref{P:ae1}. 

If two well-ordered spectra ${\mathcal S}_{G} = \{ G_{\alpha}, p_{\alpha}^{\alpha +1},\tau \}$ and ${\mathcal S}_{T} = \{ T_{\alpha}, q_{\alpha}^{\alpha +1},\tau \}$ are given and we wish to investigate properties of a map $f \colon \lim {\mathcal S}_{G} \to \lim {\mathcal S}_{T}$, even an ability to obtain a complete information hidden in spaces {\em and} short projections of these spectra will not suffice. This is not surprising at all since spectra ${\mathcal S}_{G}$ and ${\mathcal S}_{T}$, while describing domain and range of $f$, have nothing to do with $f$ itself. After accepting this conclusion one would naturally make the next step which would basically consist of representing the map $f$ itself as the limit of ``level maps". In other words we need to find a morphism 
\[ \{ f_{\alpha} \colon G_{\alpha} \to T_{\alpha}, \tau \} \colon {\mathcal S}_{G} \to {\mathcal S}_{T} \]
limit of which is $f$. Here is the corresponding diagram

\bigskip

\[ \begin{CD}
\lim {\mathcal S}_{G} @>>> \cdots @>>> G_{\alpha +1} @>p_{\alpha}^{\alpha +1}>> G_{\alpha} @>>> \cdots @>>> G_{0}\\
@V{f}VV @VVV @V{f_{\alpha +1}}VV @VV{f_{\alpha}}V  @VVV @VV{f_{0}}V\\
\lim {\mathcal S}_{T} @>>> \cdots @>>> T_{\alpha +1} @>q_{\alpha}^{\alpha +1}>> T_{\alpha} @>>> \cdots @>>> T_{0}
\end{CD}\]
\begin{center}
\textsc{Diagram 1} 
\end{center}

Representation of $f$ as the limit of such a morphism is not always possible. Trivial counterexamples can be found for maps between {\em metrizable} compacta. But according to ``the effect of uncountability" discovered by \v{S}\v{c}epin, any map between {\em nonmetrizable} compacta admits such a representation (this is basically what the Spectral Theorem states). Having now information about $f_{\alpha}$'s available in principle we find ourselves in a situation described above. This information will not suffice to restore properties of $f$. This is information about {\em objects} and information needed is, as we have already seen, hidden in {\em morphisms}. What are objects and what are morphisms in this case? Clearly $f_{\alpha}$'s are objects (representing the given object $f$) and morphisms are nothing else but commutative square diagrams

\[
\begin{CD}
G_{\alpha +1} @>f_{\alpha +1}>> T_{\alpha +1}\\
@V{p_{\alpha}^{\alpha +1}}VV  @VV{q_{\alpha}^{\alpha +1}}V \\
G_{\alpha} @>f_{\alpha}>> T_{\alpha} .
\end{CD}
\]
\begin{center}
\textsc{Diagram 2}
\end{center}

\noindent The rest is a matter of technicalities. How one would investigate properties of such square diagrams? A possible approach is through their {\em characteristic} maps. In the above notations this is the diagonal product $\chi \colon p_{\alpha}^{\alpha +1} \triangle f_{\alpha +1} \colon G_{\alpha +1} \to K$, where $K$ is the pullback

\[ K = \{ (g,t) \in G_{\alpha} \times T_{\alpha +1} \colon f_{\alpha}(g) = q_{\alpha}^{\alpha +1}(t) \} \]

\noindent of maps $f_{\alpha}$ and $q_{\alpha}^{\alpha +1}$. Here is the corresponding diagram
\bigskip

\begin{picture}(280,150)
\put(96,0){$G_{\alpha}$}
\put(90,130){$G_{\alpha +1}$}
\put(235,0){$T_{\alpha}$}
\put(230,130){$T_{\alpha +1}$}
\put(119,133){\vector(1,0){105}}
\put(102,125){\vector(0,-1){110}}
\put(240,125){\vector(0,-1){110}}
\put(170,139){$f_{\alpha +1}$}
\put(175,8){$f_{\alpha}$}
\put(77,65){$p_{\alpha}^{\alpha +1}$}
\put(245,65){$q_{\alpha}^{\alpha +1}$}
\put(110,126){\vector(1,-1){50}}
\put(160,60){\vector(-1,-1){50}}
\put(177,77){\vector(1,1){50}}
\put(119,2){\vector(1,0){110}}
\put(163,65){$K$}
\put(140,101){$\chi$}
\put(140,31){$\operatorname{pr}_{1}|K$}
\put(177,110){$\operatorname{pr}_{2}|K$}
\end{picture}

\begin{center}
\textsc{Diagram 3}
\end{center}

\noindent Obviously $\chi$ measures how close the above diagram is to the pullback square (or the Cartesian square in the terminology of \cite{book}). 

It is easy to see that the surjectivness of all vertical arrows in Diagram 1 does not imply surjectivness of $f$. Explanation is simple: these are requirements about $f_{\alpha}$'s, i.e. {\em objects}, but not about mor\-p\-hisms. Having this in mind, it is clear that requesting additionally surjectivness of characteristic maps of all square diagrams in Diagram 1 (thus imposing conditions on {\em morphisms}), one would now expect that the limit would also be surjective. This is indeed true. Moreover it can be shown that this is the only right way to proceed: if one starts with surjective $f$, then ``almost all" square diagrams formed by {\em any} morphism representing $f$ would have surjective characteristic maps. Same applies to the openness of $f$. Requesting openness of not only $f_{\alpha}$'s but of all characteristic maps of all square diagrams we arrive to the $0$-softness of the limit map $f$. Conversely, $0$-soft maps are precisely those which can be obtained in this way. Noting that there are open but not $0$-soft maps, we see where the difference between Theorems A and E comes from. This also shows how Lemma \ref{L:parametricexp} is related to Proposition \ref{P:exp}.

Definitions of few remaining concepts are discussed below. For the readers convenience we mainly cite \cite{book} where a complete discussion, original sources and omitted details can be found.  

\section{Characterization of torus groups}\label{S:results}
The following simple and well-known fact plays an important role below. A standard proof is based on Pontryagin-van Kampen duality (I thank K.~H.~Hofmann for pointing this to me out). Here we present its topological translation.

\begin{lem}\label{L:parametric}
Let $p \colon G \to T$ be a continuous surjective homomorphism between compact Abelian groups. If $\ker (p)$ is a metrizable $AE(1)$-com\-pac\-tum, then $p$ is (in both topological and algebraic senses) a trivial fibration.
\end{lem}
\begin{proof}
Since the Abelian group $\ker (p)$ is a metrizable $AE(1)$-com\-pac\-tum, it is a torus group (the equivalence $(11) \Longleftrightarrow (12)$ of Theorem A). Since torus groups are injective objects in the category of compact Abelian groups (the equivalence $(10) \Longleftrightarrow (12)$ of Theorem A), there exists a continuous retraction $r \colon G \to \ker (p)$ which is a homomorphism. Obviously, there is an isomorphism $h \colon G \to \ker (r) \times \ker (p)$ (here is the explicit formula: $h(g) = (g - r(g), r(g))$, $g \in G$). Now observe that $q = p|\ker (r) \colon \ker (r) \to T$ is an isomorphism. Consequently, the following diagram

\[
\begin{CD}
G @>h>> \ker (r) \times \ker (p)\\
@V{p}VV @VV{\operatorname{pr}_{1}}V\\
T @>q^{-1}>> \ker (r)
\end{CD}
\]

\noindent commutes. This proves the lemma.
\end{proof}

The following statement provides a characterization of $AE(1)$-com\-pact Abelian groups in terms of well-ordered inverse spectra. It should be especially emphasized that characterization of these groups in terms of $\omega$-spectra is an immediate consequence of \cite[Proposition 6.3.5 and Lemma 8.2.1]{book} and \v{S}\v{c}epin's Spectral Theorem for compacta \cite[Theorem 1.3.4]{book}. But for establishing our results we need well-ordered spectra - they allow us to proceed by induction. On the other hand for well ordered spectra we are forced to reprove already known topological versions in the presence of group structure. The latter only slightly simplifies a typical proof (see \cite{book}) although the outcome is much stronger.

\begin{pro}\label{P:ae1}
Let $G$ be a compact Abelian group of weight $\tau > \omega$. Then the following conditions are equivalent:
\begin{itemize}
\item[(a)]
$G$ is the torus $\text{\f T}^{\tau}$.
\item[(b)]
$G$ is an $AE(1)$-compactum.
\item[(c)]
There exists a well-ordered inverse spectrum ${\mathcal S}_{G} = \{ G_{\alpha}, p_{\alpha}^{\alpha +1},\tau \}$ satisfying the following properties:
\begin{enumerate}
\item
$G_{\alpha}$ is a compact Abelian group and $p_{\alpha +1} \colon G_{\alpha +1} \to G_{\alpha}$ is a continuous homomorphism, $\alpha < \tau$.
\item
If $\beta < \tau$ is a limit ordinal, then the diagonal product 
\[\triangle\{ p_{\alpha}^{\beta} \colon \alpha < \beta\} \colon G_{\beta} \to \lim\{ G_{\alpha}, p_{\alpha}^{\alpha+1}, \alpha < \beta \}\]
is an isomorphism.
\item
$G$ is isomorphic to $\lim {\mathcal S}$.
\item
The short projection $p_{\alpha}^{\alpha +1} \colon G_{\alpha +1} \to G_{\alpha}$ is isomorphic to the trivial fibration $G_{\alpha} \times \ker\left(p_{\alpha}^{\alpha +1}\right) \to G_{\alpha}$ the fiber $\ker\left(p_{\alpha}^{\alpha +1}\right)$ of which is a metrizable torus, $\alpha < \tau$.
\item
$G_{0}$ is a metrizable torus. 
\end{enumerate}
\end{itemize}
\end{pro}
\begin{proof}
$(a) \Longrightarrow (b)$. Note that $\text{\f T}^{\tau}$ is an $AE(1)$-compactum.

$(c) \Longrightarrow (a)$. Apply a straightforward transfinite induction to conclude that $G$ is topologically and algebraically equivalent to the product $G_{0} \times \{ \ker\left( p_{\alpha}^{\alpha +1}\right) \colon \alpha < \tau \}$. Properties 4 and 5 show that such a product is isomorphic to the torus group $\text{\f T}^{\tau}$.

$(b) \Longrightarrow (c)$. We may assume that $G$ is a closed subgroup of a product $\prod\{ X_{a} \colon a \in A\}$, $|A| = \tau$, of compact metrizable Abelian groups $X_{a}$, $a \in A$, each of which is connected and locally connected (clearly each of this $X_{a}$'s can be chosen to be the circle group $\text{\f T}$ but in order to avoid a confusion we deliberately use a different notation). By \cite[Theorem 4.1]{dran84}, there exists an open $1$-invertible map $f \colon Y \to \prod\{ X_{a} \colon a \in A \}$, where $Y$ is a $1$-dimensional compactum of weight $\tau$. Consider the inverse image $f^{-1}(G) \subseteq Y$ of $G$ and the map $f|f^{-1} \colon f^{-1}(G) \to G$. Since $\dim Y = 1$ and since $G$, according to (b), is an $AE(1)$-compactum, there exists a map $g \colon Y \to G$ such that $g|f^{-1}(G) = f|f^{-1}(G)$.

Next let us denote by 
\[ \pi_{B} \colon \prod\{ X_{a} \colon a \in A\} \to \prod\{ X_{a} \colon a \in B\}\] 
and
\[\pi^{B}_{C} \colon \prod\{ X_{a} \colon a \in B\} \to \prod\{ X_{a} \colon a \in C \}\]
the natural projections onto the corresponding subproducts ($C \subseteq B \subseteq A$).
We call a subset $B \subseteq A$ {\em admissible} (compare with the proof of \cite[Theorem 6.3.1]{book}) if the following equality
\[ \pi_{B}(g(f^{-1}(x))) = \pi_{B}(x)\]
is true for each point $x \in \pi_{B}^{-1}\left(\pi_{B}(G)\right)$. We need the following properties of admissible sets.

{\em Claim 1. The union of arbitrary collection of admissible sets is admissible}. 

Indeed let $\{ B_{t} \colon t \in T\}$ be a collection of admissible sets and $B = \cup \{ B_{t} \colon t \in T\}$. Let $x \in \pi_{B}^{-1}\left(\pi_{B}(G)\right)$. Clearly $x \in \pi_{B_{t}}^{-1}\left(\pi_{B_{t}}(G)\right)$ for each $t \in T$ and consequently 
\[ \pi_{B_{t}}(g(f^{-1}(x))) = \pi_{B_{t}}(x) \;\;\text{for each}\;\; t \in T .\]
Obviously, $\pi_{B}(x) \in \pi_{B}(g(f^{-1}(x)))$ and it therefore suffices to show that the set $\pi_{B}(g(f^{-1}(x)))$ contains only one point. Assuming that there is a point $y \in \pi_{B}(g(f^{-1}(x)))$ such that $y \neq \pi_{B}(x)$ we conclude (having in mind that $B = \cup\{ B_{t} \colon t \in T\}$) that there must be an index $t \in T$ such that $\pi_{B_{t}}^{B}(y) \neq \pi_{B_{t}}^{B}\left(\pi_{B}(x)\right)$. But this is impossible
\[ \pi_{B_{t}}^{B}(y) \in \pi_{B_{t}}^{B}\left(\pi_{B}(g(f^{-1}(x)))\right) = \pi_{B_{t}}(g(f^{-1}(x))) = \pi_{B_{t}}(x) = \pi_{B_{t}}^{B}\left(\pi_{B}(x)\right) .\]

{\em Claim 2. If $B \subseteq A$ is admissible, then the restriction $\pi_{B}|G \colon G \to \pi_{B}(G)$ is $1$-soft}.

Let $\varphi \colon Z \to \pi_{B}(G)$ and $\varphi_{0} \colon Z_{0} \to G$ be two maps defined on a $1$-dimensional compactum $Z$ and its closed subset $Z_{0}$ respectively. Assume that $\pi_{B}\varphi_{0} = \varphi |Z_{0}$. We wish to construct a map $\phi \colon Z \to G$ such that $\phi |Z_{0} = \varphi_{0}$ and $\pi_{B}\phi = \varphi$, i.e. $\phi$ makes the diagram

\begin{picture}(300,150)
\put(80,110){$Z$}
\put(80,10){$Z_{0}$}
\put(80,100){$\cap$}
\put(87,102){\vector(0,-1){79}}
\put(65,60){$\operatorname{incl}$}
\put(95,17){\vector(1,1){90}}
\put(95,112){\vector(1,0){90}}
\put(95,12){\vector(1,0){80}}
\put(190,110){$G$}
\put(180,10){$\pi_{B}(G)$}
\put(195,107){\vector(0,-1){83}}
\put(135,118){$\varphi_{0}$}
\put(135,17){$\varphi$}
\put(135,67){$\phi$}
\put(198,60){$\pi_{B}|G$}
\end{picture}

\noindent commutative. Since, according to our choice, all $X_{a}$'s are $AE(1)$-compacta, so is the product $\prod\{ X_{a} \colon a \in A-B\}$. This implies the $1$-softness of the projection $\pi_{B}$ and hence of its restriction 
\[ \pi_{B}|\pi_{B}^{-1}\left(\pi_{B}(G)\right) \colon \pi_{B}^{-1}\left(\pi_{B}(G)\right) \to \pi_{B}(G) .\] 
Then there exists a map $\phi^{\prime\prime} \colon Z \to \pi_{B}^{-1}\left(\pi_{B}(G)\right)$ such that $\phi^{\prime\prime}|Z_{0} = \varphi_{0}$ and $\pi_{B}\phi^{\prime\prime} = \varphi$. Since $f$ is $1$-invertible (and $\dim Z = 1$), there exists a map $\phi^{\prime} \colon Z \to Y$ such that $f\phi^{\prime} = \phi^{\prime\prime}$. Now let $\phi = g\phi^{\prime}$. Since $g|f^{-1}(G) = f|f^{-1}(G)$, we have $\varphi_{0} = \phi |Z_{0}$. Finally observe that the admissibility of $B$ implies $\varphi = \pi_{B}\phi$ as required.

{\em Claim 3. For each countable subset $C \subseteq A$ there exists a countable admissible subset $B \subseteq A$ such that $C \subseteq B$}.

Since $w(Y) = \tau$, it follows (consult \cite{book}) that $Y$ can be represented as the limit space of an $\omega$-spectrum ${\mathcal S}_{Y} = \{ Y_{B}, q^{B}_{C}, \exp_{\omega}A \}$ consisting of metrizable compacta $Y_{B}$, $B \in \exp_{\omega}A$, and continuous surjections $q_{C}^{B} \colon Y_{B} \to Y_{C}$, $C \subseteq B$, $C,B \in \exp_{\omega}A$. Consider also the standard $\omega$-spectrum ${\mathcal S}_{X} = \{ \prod\{ X_{a} \colon a \in B\} ,\pi_{C}^{B}, \exp_{\omega}A \}$ consisting of countable subproducts of the product $\prod\{ X_{a} \colon a \in A \}$ and corresponding natural projections. Obviously the full product coincides with the limit of ${\mathcal S}_{X}$. One more $\omega$-spectrum arises naturally. This is the spectrum ${\mathcal S}_{G} = \{ \pi_{B}(G), \pi_{C}^{B}|\pi_{B}(G), \exp_{\omega}A \}$ the limit of which coincides with $G$. 

Consider the map $f \colon \lim {\mathcal S}_{Y} \to \lim {\mathcal S}_{X}$. Applying the Spectral Theorem (\cite[Theorem 1..3.4]{book}) there is a cofinal and $\omega$-closed subset ${\mathcal T}_{1}$ of $\exp_{\omega}A$ such that for each $B \in {\mathcal T}_{f}$ there is a map $f_{B} \colon Y_{B} \to \prod\{ X_{a} \colon a \in B \}$ such that $f_{B}q_{B} = \pi_{B}f$. Moreover, these maps form a morphism
\[ \{ f_{B} ; B \in {\mathcal T}_{f}\} \colon {\mathcal S}_{Y} \to {\mathcal S}_{X}\]
limit of which coincides with $f$. Since $f$ is open (and closed), we may assume without loss of generality (considering a smaller cofinal and $\omega$-subset of ${\mathcal T}_{f}$ if necessary) that the above indicated morphism is bicommutative. This simply means that $q_{B}f^{-1}(K) = f_{B}^{-1}\left(\pi_{B}(K)\right)$ for any $B \in {\mathcal T}_{f}$ and any closed subset $K$ of the product $\prod\{ X_{a} \colon a \in A\}$.

Similarly applying the Spectral Theorem to the map $g \colon \lim {\mathcal S}_{Y} \to \lim {\mathcal S}_{G}$ we obtain a cofinal and $\omega$-closed subset ${\mathcal T}_{g}$ of $\exp_{\omega}A$ and the associated to it morphism
\[ \{ g_{B} \colon Y_{B} \to \pi_{B}(G) ; B \in {\mathcal T}_{g} \} \colon {\mathcal S}_{Y} \to {\mathcal S}_{G}\]
limit of which coincides with the map $g$.

By \cite[Proposition 1.1.27]{book}, the intersection ${\mathcal T} = {\mathcal T}_{f} \cap {\mathcal T}_{g}$ is still a cofinal and $\omega$-closed subset of $\exp_{\omega}A$. It therefore suffices to show that each $B \in {\mathcal T}$ is an admissible subset of $A$. Consider a point $x \in \pi_{B}^{-1}\left( \pi_{B}(G)\right)$. First observe that bicommutativity of the morphism associated with ${\mathcal T}_{f}$ implies that $q_{B}(f^{-1}(x)) = f_{B}^{-1}(\pi_{B}(x))$. Since the maps $f_{B}$ and $g_{B}$ coincide on $f_{B}^{-1}(\pi_{B}(G))$ we have 
\begin{multline*}
\pi_{B}(g(f^{-1}(x))) = g_{B}(q_{B}(f^{-1}(x))) = g_{B}(f_{B}^{-1}(\pi_{B}(x))) =\\ f_{B}(f_{B}^{-1}(\pi_{B}(x))) = \pi_{B}(x)
\end{multline*}
as required.

{\em Claim 4. If $C$ and $B$ are admissible subsets of $A$ and $C \subseteq B$, then the map $\pi_{C}^{B}|\pi_{B}(G) \colon \pi_{B}(G) \to \pi_{C}(G)$ is $1$-soft}.

This property follows from Claim 2 and \cite[Lemma 6.1.15]{book}.

After having all the needed properties of admissible subsets established we proceed as follows. Since $|A| = \tau$ we can write $A = \{ a_{\alpha} \colon \alpha < \tau\}$. By Claim 3, each $a_{\alpha} \in A$ is contained in a countable admissible subset $B_{\alpha} \subseteq A$. Let $A_{\alpha} = \cup\{ B_{\beta} \colon \beta \leq \alpha\}$. We use these sets to define a transfinite inverse spectrum ${\mathcal S} = \{ G_{\alpha}, p_{\alpha}^{\alpha +1}, \tau \}$ as follows. Let $G_{\alpha} = \pi_{A_{\alpha}}(G)$ and $p_{\alpha}^{\alpha +1} = \pi_{A_{\alpha}}^{A_{\alpha +1}}|G_{\alpha +1}$ for each $\alpha < \tau$. Properties 1-3 of the spectrum ${\mathcal S}_{G}$ are satisfied by construction. Since, $G_{0}$ is a metrizable $AE(1)$-compactum, it follows from the equivalence $(11) \Longleftrightarrow (12)$ of Theorem A, that $G_{0}$ is a metrizable torus. Finally, property 4 is a consequence of Claim 4 and Lemma \ref{L:parametric}. Proof is completed.
\end{proof}

Next we investigate the map $\exp_{G} \colon L(G) \to G$.

\begin{lem}\label{L:parametricexp}
Let $p \colon G \to T$ be a continuous homomorphism of compact Abelian groups and $\ker (p)$ is a metrizable $AE(1)$-compactum. Then the diagram
\[
\begin{CD}
L(G) @>\exp_{G}>> G\\
@V{L(p)}VV  @VV{p}V\\
L(T) @>\exp_{T}>> T
\end{CD}
\]

\noindent is $0$-soft.
\end{lem}
\begin{proof}
By Lemma \ref{L:parametric}, we may without loss of generality assume that $G$ is the product $T \times \ker (p)$ and the homomorphism $p$ coincides with the projection $\pi_{T} \colon T \times \ker (p) \to T$. Now consider the following commutative diagram (compare it with Diagram 3):

\begin{picture}(280,150)
\put(90,0){$L(T)$}
\put(90,130){$L(G)$}
\put(255,0){$T$}
\put(240,130){$T \times \ker (p) = G$}
\put(119,133){\vector(1,0){105}}
\put(105,127){\vector(0,-1){110}}
\put(260,127){\vector(0,-1){110}}
\put(175,139){$\exp_{G}$}
\put(175,8){$\exp_{T}$}
\put(80,65){$L(p)$}
\put(265,65){$\pi_{T} = p$}
\put(110,126){\vector(1,-1){45}}
\put(165,65){\vector(-1,-1){50}}
\put(205,85){\vector(1,1){40}}
\put(119,2){\vector(1,0){130}}
\put(153,70){$L(T) \times \ker (p)$}
\put(140,101){$\chi$}
\put(140,33){$\pi_{L(T)}$}
\put(180,110){$\exp_{T} \times \operatorname{id}$}
\end{picture}

\bigskip

\noindent in which $\chi = L(p) \triangle \exp_{G}$ is the characteristic map of the previous diagram \cite[Definition 6.2.1]{book}. We therefore need to show the $0$-softness of $\chi$. A straightforward verification shows that in this situation $L(G)$ is isomorphic to the product $L(T) \times L(\ker (p))$ and the homomorphism $\chi$ coincides with the product $\operatorname{id} \times \exp_{\ker (p)}$. Observe also that with these identifications $L(p)$ becomes isomorphic to the projection $L(T) \times L(\ker (p)) \to L(T)$. Since $\ker (p)$ is a metrizable $AE(1)$-compactum, the map $\exp_{\ker (p)}$ is an open surjection between completely metrizable separable spaces (equivalence $(3) \Longleftrightarrow (4) \Longleftrightarrow (11)$). By \cite[Corollary 6.1.27]{book}, $\exp_{G}$ is $0$-soft. This implies $0$-softness of $\chi$.
\end{proof}

\begin{pro}\label{P:exp}
Let $G$ be a compact connected Abelian group. Then the following conditions are equivalent:
\begin{itemize}
\item[(a)]
$G$ is an $AE(1)$-compactum.
\item[(b)]
The map $\exp_{G} \colon L(G) \to G$ is $0$-soft.
\end{itemize}
\end{pro}
\begin{proof}
Since the statement is known to be true in the metrizable case we assume that $w(G) = \tau > \omega$. 

$(a) \Longrightarrow (b)$. Represent $G$ as the limit of a well-ordered inverse spectrum ${\mathcal S} = \{ G_{\alpha}, p_{\alpha}^{\alpha +1},\tau\}$ with the properties listed in Proposition \ref{P:ae1}(c). The space $L(G)$ is then the limit of the associated spectrum $L({\mathcal S}) = \{ L(G_{\alpha}), L\left( p_{\alpha}^{\alpha +1}\right) ,\tau \}$ and the map $\exp_{G} \colon L(G) \to G$ is the limit of the morphism
\[ \{ \exp_{G_{\alpha}} \colon L(G_{\alpha}) \to G_{\alpha}, \tau\} \colon L({\mathcal S}) \to {\mathcal S} .\]
According to Lemma \ref{L:parametricexp} this morphism is $0$-soft. $0$-softness of the limit map $\exp_{G}$ follows now from \cite[Theorem 6.3.1]{book} (applied with $n=0$).

$(b) \Longrightarrow (a)$. The proof of this part follows the proof of \cite[Theorem 6.3.1]{book}. Only two things have to be taken into account. The first has already been mentioned: the space $L(G)$ is isomorphic to $\text{\f R}^{\operatorname{rank}(\widehat{G})}$ and hence is an $AE(0)$-space in the sense of \cite[Definition 6.1.12]{book}. The second: $G$ being a compact group is also an $AE(0)$-compactum (Theorem C). Following the indicated proof and making slight adjustments (as in the proof of Proposition \ref{P:ae1}) we arrive to the following situation. There exists a well-ordered spectrum ${\mathcal S}_{G} = \{ G_{\alpha}, p_{\alpha}^{\alpha +1}, \tau\}$, satisfying properties 1-3 and 5 of Proposition \ref{P:ae1}(c) and a part of property 4 which guarantees that the fiber $\ker\left( p_{\alpha}^{\alpha +1}\right)$ is metrizable. Also there exists a morphism (we keep notation of the first part of the proof)
\[ \{ \exp_{G_{\alpha}} \colon L(G_{\alpha}) \to G_{\alpha}, \tau\} \colon L({\mathcal S}) \to {\mathcal S} \]
between this spectra which is $0$-soft in the sense that characteristic maps of all the square diagrams

\[
\begin{CD}
L(G_{\alpha +1}) @>\exp_{G_{\alpha +1}}>> G_{\alpha +1}\\
@V{L(p_{\alpha}^{\alpha +1})}VV  @VV{p_{\alpha}^{\alpha +1}}V\\
L(G_{\alpha}) @>\exp_{G_{\alpha}}>> G_{\alpha}
\end{CD}
\]
are $0$-soft. In addition, it follows from the above remarks that typical projection $L(p_{\alpha}^{\alpha +1})$ is isomorphic to the projection $\operatorname{pr}_{1} \colon L(G_{\alpha}) \times \text{\f R}^{\omega} \to L(G_{\alpha})$. Since the characteristic map of the above diagram is $0$-soft, it is surjective \cite[Lemma 6.1.13]{book}. This means that $\left( p_{\alpha}^{\alpha +1}\right)^{-1}(\exp_{G_{\alpha}}(x)) = \exp_{G_{\alpha +1}}\left( L^{-1}(p_{\alpha}^{\alpha +1})(x)\right)$ for any point $x \in L(G_{\alpha})$. Apply the latter to the point ${\mathbf 0} \in L(G_{\alpha})$. We have $\ker\left( p_{\alpha}^{\alpha +1}\right) = \exp_{G_{\alpha +1}}\left( L^{-1}(p_{\alpha}^{\alpha +1})({\mathbf 0})\right)$. Since, as was noted, the map $L(p_{\alpha}^{\alpha +1})$ is isomorphic to the projection $\operatorname{pr}_{1}$, it follows that $L^{-1}(p_{\alpha}^{\alpha +1})({\mathbf 0})$ is isomorphic to $\text{\f R}^{\omega}$. Consequently a metrizable compact Abelian group $\ker\left( p_{\alpha}^{\alpha +1}\right)$ being a {\em surjective} image of $\text{\f R}^{\omega}$ (through the map $\exp_{G_{\alpha +1}}$) is arcwise connected. Equivalence $(1) \Longleftrightarrow (12)$ applied to $\ker\left( p_{\alpha}^{\alpha +1}\right)$ shows that the latter is a metrizable torus. This finishes verification of property 4 of Proposition \ref{P:ae1}(c). Since $G_{0}$ is obviously a metrizable torus, it follows (Proposition \ref{P:ae1}) that $G$ is an $AE(1)$-compactum. Proof is completed.
\end{proof}

Theorem E from the introduction is a consequence of Propositions \ref{P:ae1} and \ref{P:exp}.

\section{Free Abelian Groups}\label{S:free}
Discussion of spectral analysis in Section \ref{S:prelim} indicates that although the Spectral Theorem has originally been designed for topological category, its explicit categorical nature allows us to obtain corresponding counterparts in other categories (such as category of compact (Abelian) groups, Shape category, category of spaces admitting effective actions of compact groups etc. \cite{book}) as well. Here we consider the category of Abelian groups. 

Let $\tau \geq \omega$ be a given cardinal number. Every Abelian group $G$ can be represented as the union $\cup\{ G_{\alpha} \colon \alpha \in A \}$ of its subgroups $G_{\alpha}$, $\alpha \in A$, of cardinality $|G_{\alpha}| \leq \tau$, where the partially ordered and directed indexing set $A$ is $\tau$-complete (this basically means that $A$ is not only directed but also $\tau$-directed - contains elements majorating all elements of any subset of cardinality $\leq \tau$; see \cite{book} for details). Basic example of such an indexing set $A$ is the one ($\exp_{\tau}G$) corresponding to the representation of $G$ as the union of {\em all} subgroups of cardinality $\leq \tau$. Other $\tau$-complete indexing sets would be (at least in all the situations we deal with) cofinal and $\tau$-closed (contain supremums of subsets of cardinality $\leq \tau$) subsets of $\exp_{\tau}G$.

Let us say that the above representation is $\tau$-smooth\footnote{This terminology is suggested by \cite{eklof}.} if, in addition, $G_{\beta} = \cup\{ G_{\alpha} \colon \alpha \in B\}$, whenever $B$ is a chain in $A$, $|B| \leq \tau$ and $\beta = \sup B$. Observe that $\tau$-smooth decompositions (i.e., direct $\tau$-spectra) are in one to one correspondence with $\tau$-spectra of compact Abelian groups.

\begin{thm}[$\tau$-Spectral Theorem]\label{T:smooth}
Let $f \colon G \to L$ be a homomorphism between Abelian groups. Suppose also that $\{ G_{\alpha} \colon \alpha \in A\}$ and $\{ L_{\alpha} \colon \alpha \in A \}$ are {\em any} two $\tau$-smooth decompositions of $G$ and $L$ with the same indexing set $A$. Then there exists a cofinal and $\tau$-closed subset $B$ of $A$ such that $f(G_{\alpha}) \subseteq L_{\alpha}$ for each $\alpha \in B$.
\end{thm}
\begin{proof}
Apply the Spectral Theorem for $\tau$-spectra of compact Abelian groups (see \cite[Lemma 8.2.1]{book} where the non Abelian and non compact version is proved, consult also \cite{shch}), note that $\tau$-smooth decompositions are nothing else but $\tau$-spectra and apply duality. 
\end{proof}

There is no difficulty in defining the concept corresponding to well-ordered spectra of compact Abelian groups. These are decompositions of Abelian groups into increasing well-ordered systems of subgroups which are called {\em smooth} \cite{eklof}, i.e., $\{ G_{\alpha} \colon \alpha < \tau \}$ is a smooth decomposition of $G$ if $G_{\alpha} \subseteq G_{\alpha +1}$, $\alpha < \tau$, and if $G_{\beta} = \cup\{ G_{\alpha} \colon \alpha < \beta \}$ for any limit ordinal $\beta < \tau$. 

\begin{thm}[Eklof,\cite{eklof}]\label{ekl}
Let $\{ G_{\alpha} \colon \alpha < \tau \}$ be a smooth decomposition of an Abelian group $G$ such that $G_{0}$ is free and $G_{\alpha +1}/G_{\alpha}$ is free for every $\alpha < \tau$. Then $G$ is free.
\end{thm}
The above statement in essence gives one half of characterization of free Abelian groups. This half corresponds to the implication $(c) \Longrightarrow (a)$ of Proposition \ref{P:ae1}. Since, as we have seen, the other implication $(a) \Longrightarrow (c)$ of Proposition \ref{P:ae1} is also true, we obtain the following result.

\begin{thm}\label{T:free}
An (uncountable) Abelian group $G$, $|G| = \tau$, is free if and only if it admits a smooth decomposition $\{ G_{\alpha} \colon \alpha < \tau \}$ such that $G_{0}$ is {\em countable} and free and $G_{\alpha +1}/G_{\alpha}$ is {\em countable} and free for every $\alpha < \tau$.
\end{thm}

From methodological point of view it is perhaps of some interest to note that the well known \cite[Theorem 4,p.45]{lang} statement ``a subgroup of a free Abelian group is free" has in fact a countable nature.

\begin{cor}\label{C:subgroup}
Statements (a) and (b) are equivalent:
\begin{itemize}
\item[(a)]
A subgroup of a countable free Abelian group is free.
\item[(b)]
A subgroup of free Abelian group is free.
\end{itemize}
\end{cor}
\begin{proof}
Let $L$ be a subgroup of an Abelian group $G$. If $G$ is free it has a smooth decomposition $\{ G_{\alpha} \colon \alpha < \tau \}$ with the properties indicated in Theorem \ref{T:free}. Then $\{ L_{\alpha} = L \cap G_{\alpha} \colon \alpha < \tau \}$ is a smooth decomposition of $L$. Since $L_{0} \subseteq G_{0}$ and $L_{\alpha +1}/L_{\alpha} \subseteq G_{\alpha +1}/G_{\alpha}$, $\alpha < \tau$, we conclude (by (a)) that $\{ L_{\alpha} \colon \alpha < \tau \}$ satisfies conditions of Theorem \ref{T:free} and hence $L$ is free.
\end{proof}

\begin{cor}\label{C:image}
A continuous and homomorphic image of a torus is a torus.
\end{cor}

\end{document}